\documentclass[]{article}
\usepackage{latexsym,graphics}

\newcommand{\kitem}{\begin{itemize}\vspace{-2ex}}
\newcommand{\kenditem}{\vspace{-1ex}\end{itemize}}

\newcommand{\Q}{{I\!\!\!\!Q}}
\newcommand{\R}{{I\!\!R}}

\newcommand{\ku}{\underline}

\newcommand{\kss}{\scriptscriptstyle}

\newcommand{\kb}{{\kss \bullet}}
\newcommand{\ktype}{\sc}

\reversemarginpar
\newcommand{\kkk}[1]{}

\newcommand{\vv}{{v}}
\newcommand{\vvg}{\tilde{\vv}}
\newcommand{\vexpo}{e}
\newcommand{\ve}[1]{{\vexpo_{#1}}}
\newcommand{\veg}{\tilde{\vexpo}}
\newcommand{\vsero}{s}
\newcommand{\vs}[1]{{\vsero_{#1}}}
\newcommand{\vp}{p}
\newcommand{\vpk}[1]{{\vp({#1})}}
\newcommand{\vpi}[1]{{\vp_{#1}}}
\newcommand{\vqq}{q}
\newcommand{\vq}[1]{{\vqq_{#1}}}
\newcommand{\di}{{i}}
\newcommand{\dk}{{k}}
\newcommand{\va}{a}
\newcommand{\vA}{A}
\newcommand{\vn}{n}
\newcommand{\vak}[1]{{\va({#1})}}
\newcommand{\vai}[1]{\va_{#1}}
\newcommand{\vag}[1]{\tilde{\va}_{#1}}
\newcommand{\vc}{c}
\newcommand{\vf}{f}
\newcommand{\vg}{g}

\newcounter{Abschnitt}[section]
\newcommand{\neu}[1]{\protect\refstepcounter{Abschnitt}\protect
   \label{#1}\vspace{1ex}
   {\bf (\protect\arabic{section}.\protect\arabic{Abschnitt})}
                     $\qquad$\kkk{#1}}
\newcommand{\zitat}[2]{(\protect\ref{#1}.\protect\ref{#1-#2})}

\begin{document}
\title{Estimating Vaccine Coverage by Using Computer Algebra}
\author{Doris Altmann \quad\qquad Klaus Altmann }
\date{}
\maketitle
\begin{abstract}
The approach of N.~Gay for estimating the coverage of a multivalent vaccine
from 
antibody prevalence data in certain age cohorts
is improved by using computer aided elimination theory of variables.
Hereby, Gay's usage of numerical approximation can be replaced by 
exact formulas which are surprisingly nice, too.
\end{abstract}

%
%
\section{Introduction}\label{intro}


\neu{intro-overview} 
Nigel~Gay \cite{Gay} has estimated the coverage of MMR 
(measles, mumps, rubella) multivalent vaccination 
in a fixed age cohort by the following method:\\
The rates $\vpk{\pm,\pm,\pm}$ of being seropositive with
each of the three diseases depend, via a polynomial system $F$,  
on the MMR coverage $\vv$, the exposition factors $\ve{\di}$,
and the rates $\vs{\di}$ of seroconversion; the index $\di=1,2,3$ stands for
measles, mumps and rubella, respectively. 
On the other hand, it is the $\vpk{\pm,\pm,\pm}$ which can be
obtained from the available data. Hence, a maximum likelihood
approach provides estimations of $\vv$, $\ve{\di}$, and $\vs{\di}$.
\par

Gay's approach leads to numerical methods
of finding values $\vv,\ve{\di}, \vs{\di}$ that minimize the distance
between $F_{\pm,\pm,\pm}(\vv,\ve{\di}, \vs{\di})$ and the measured
$\vpk{\pm,\pm,\pm}$. The present paper replaces this part by providing
{\em exact} formulas describing the inverse of the polynomial map 
$\,F:\R^7\to\R^8$. 
Note that the image of $F$ is contained in the hyperplane
$[\sum \vpk{\pm,\pm,\pm}=1]$, i.e.\ it is 7-dimensional like the source
space of $F$.\\
The final result providing our estimation of $\vv,\ve{\di}, \vs{\di}$
may be found in Theorem \zitat{app}{theorem}.
\par


\neu{intro-diff} 
We make the same three assumptions used by Gay \cite{Gay}:
\kitem
\item[(1)]
Vaccinated children who do not seroconvert as a result of vaccination have
the same probability of being seropositive as an unvaccinated child of the
same age (i.e., $\ve{\di}$).
\item[(2)]
In a single individual, seroconversion to each vaccine component is independent.
\item[(3)]
Risk of exposure to infection is homogeneous within each age cohort and
infection with each disease is independent.
\kenditem
However, we eliminate another assumption which is silently made
in \cite{Gay} in that we do not assume that the seroconversion $\vs{\di}$ 
for the
$\di$-th disease is independent of age.
\par


\neu{intro-ack}
We would like to thank Duco van Straten for the useful discussions
concerning the exciting mathematical pattern hidden in the
MMR problem and its solution. Moreover, we are greatful to Nigel Gay
for sending us his manuscript 
including the data of the ESEN (European Seroepidemiological Network) Project.
\par

%
%
\section{The MMR system}\label{mmr}


\neu{mmr-recall}
First, let us recall from \cite{Gay} the involved variables and their mutual
relationship. Fixing one of the age cohorts, we denote by
\kitem
\item
$\vv$ the proportion of children who have received the multivalent vaccine
(``MMR coverage''),
\item
$\ve{\di}$ the rate measuring the exposure to natural infection
with disease $\di$
%
%
(``exposition factor''),
\item
$\vs{\di}$ the proportion of children previously with no detectable
antibody to disease $\di$ who acquire detectable antibody to disease $\di$ 
when vaccinated (``seroconversion'').
%
%
\kenditem
The rate $\vq{\di}$ measuring the presence of antibodies to disease $\di$
under the condition of being vaccinated may be easily expressed as
\[
\vq{\di} = \ve{\di}+(1-\ve{\di})\,\vs{\di}
\hspace{1em}\mbox{with } \;\di=1,2,3.
\]
From these data it is possible to obtain information about the expected
antibody prevalence in general. 
It is encoded in the $8$ variables $\vpk{\pm,\pm,\pm}$
with ``$+$'' at the $\di$-th place standing for the presence and ``$-$'' 
for the absence of antibodies to
the $\di$-th disease. Likewise, we may think about the sign triples as
numbers between $0$ (meaning ``$---$'') and $7$ (meaning ``$+++$'');
this allows the shorter description 
$\vpk{\pm,\pm,\pm}=\vpk{\dk}=\vpi{\dk}$. 
The equations are
\[
\renewcommand{\arraystretch}{1.3}
\begin{array}{rcccrcl}
\vpi{7} &=& \vpk{+,+,+} &=& \vv\,\vq{1}\,\vq{2}\,\vq{3} 
&+& (1-\vv)\,\ve{1}\,\ve{2}\,\ve{3}\\
\vpi{6} &=& \vpk{+,+,-} &=& \vv\,\vq{1}\,\vq{2}\,(1-\vq{3}) 
&+& (1-\vv)\,\ve{1}\,\ve{2}\,(1-\ve{3})\\
\vpi{5} &=& \vpk{+,-,+} &=& \vv\,\vq{1}\,(1-\vq{2})\,\vq{3} 
&+& (1-\vv)\,\ve{1}\,(1-\ve{2})\,\ve{3}\\
\vpi{4} &=& \vpk{+,-,-} &=& \vv\,\vq{1}\,(1-\vq{2})\,(1-\vq{3}) 
&+& (1-\vv)\,\ve{1}\,(1-\ve{2})\,(1-\ve{3})\\
\vpi{3} &=& \vpk{-,+,+} &=& \vv\,(1-\vq{1})\,\vq{2}\,\vq{3} 
&+& (1-\vv)\,(1-\ve{1})\,\ve{2}\,\ve{3}\\
\vpi{2} &=& \vpk{-,+,-} &=& \vv\,(1-\vq{1})\,\vq{2}\,(1-\vq{3}) 
&+& (1-\vv)\,(1-\ve{1})\,\ve{2}\,(1-\ve{3})\\
\vpi{1} &=& \vpk{-,-,+} &=& \vv\,(1-\vq{1})\,(1-\vq{2})\,\vq{3} 
&+& (1-\vv)\,(1-\ve{1})\,(1-\ve{2})\,\ve{3}\\
\vpi{0} &=& \vpk{-,-,-} &=& \vv\,(1-\vq{1})\,(1-\vq{2})\,(1-\vq{3}) 
&+& (1-\vv)\,(1-\ve{1})\,(1-\ve{2})\,(1-\ve{3})\,.
\end{array}
\]
\par

{\bf Remark:}
In \cite{Gay}, the variables $\vv$, $\ve{\di}$, $\vq{\di}$,
and $\vpi{\dk}$ carry a second index
pointing to the special age cohort; $\vs{\di}$ does not because of Gay's
assumption mentioned at the end of \zitat{intro}{diff}.
\par


\neu{mmr-prog}
The previous equations express the variables $\vpi{\dk}$ in terms of
$\vv, \ve{\di}, \vq{\di}$ or, 
since $\vs{\di}=(\vq{\di}-\ve{\di})/(1-\ve{\di})$,
in terms of $\vv, \ve{\di}, \vs{\di}$. 
Our goal is to describe the inverse dependencies, and we proceed in two steps:\\
First, using elimination theory, we produce in \zitat{mmr}{elimv} and
 \zitat{mmr}{elimrest} for each of the variables
$\vv, \ve{\di}, \vq{\di}$ a separate equation with coefficients
in the polynomial ring $\,\Q[\vpi{0},\dots,\vpi{7}]$. The surprising fact 
will be that all these equations are quadratic ones.
Then, as a second step, we will check in \zitat{mmr}{check}
which of the $2^7$ combinations actually provide a solution to our system.
The results of these investigations are gathered in 
Theorem \zitat{mmr}{check}.
\par

Before we start this program, we would like to introduce an easy technical 
trick in which we replace the variables $\vpi{\dk}$ by symbolic fractions
$\vai{\dk}/\vn$. By doing so, it changes the above equations in the obvious way.
For instance, the first one becomes
\[
\vai{7} \;=\; \vak{+,+,+} \;=\; \vn\,\vv\,\vq{1}\,\vq{2}\,\vq{3}
\;+\; \vn\,(1-\vv)\,\ve{1}\,\ve{2}\,\ve{3}\,.
\]
Since this manipulation increases both the degree and the number of
variables, it seemingly complicates the problem. 
However, using computer algebra systems, the computational time decreases
substantially.
Moreover, another advantage of our approach is that 
$\sum_{\dk=0}^7\vpi{\dk}=1$ translates into $\sum_{\dk=0}^7\vai{\dk}=\vn$.
In particular, when finally applying our formulas, we may directly
substitute the number of observed probands in each category for the
corresponding variables $\vai{\dk}$. 
The number $\vn$ equals the size of the cohort.
\par


\neu{mmr-elimv}
Let us start with eliminating $\vn,\ve{\di},\vq{\di}$ to obtain an equation for
the variable $\vv$ which is, by the way, of major interest. 
We work with the computer algebra system {\ktype Singular} 
developed at the University Kaiserslautern, \cite{Singular}.
\par

Let $R$ be a polynomial ring of characteristic zero with 16 variables 
$\vai{\dk}, \vn, \vv, \ve{\di}, \vq{\di}$.
For the monomial order we have to choose a global one, e.g.\ {\tt dp(16)}.
Transforming the 8 equations into an ideal $I\subseteq R$, the command
``{\tt eliminate(I,\vn*\vexpo(1)*\vexpo(2)*\vexpo(3)*\vqq(1)*\vqq(2)*\vqq(3))}'' 
produces a quadratic equation
\[
\vc_1(\vai{0},\dots,\vai{7}) \,\vv^2
- \vc_1(\vai{0},\dots,\vai{7}) \,\vv
+ \vc_0(\vai{0},\dots,\vai{7})\;=\;0
\]
with huge polynomials $\vc_1,\vc_0$ of degree 6 in the variables
$\vai{0},\dots,\vai{7}$.\\
We may also use {\ktype Singular} for the factorization of polynomials.
Applied to the coefficient $\vc_1$ as well as to the discriminant 
of our quadratic polynomial, this yields nice results. With
\[
\renewcommand{\arraystretch}{1.3}
\begin{array}{r@{\hspace{0.4em}}c@{\hspace{0.4em}}l}
%
\vf_1&:=&\vn \hspace{0.4em}=\hspace{0.4em}
\Big(\big(\vai{0}+\vai{3}+\vai{5}+\vai{6}\big)+
\big(\vai{7}+\vai{4}+\vai{2}+\vai{1}\big)\Big)
\vspace{1ex}\\
%
%
%
\vf_3 &:=&
\Big(\big(\vai{0}+\vai{3}+\vai{5}+\vai{6}\big)-
\big(\vai{7}+\vai{4}+\vai{2}+\vai{1}\big)\Big)
\Big(\vai{0}\vai{7}+\vai{3}\vai{4}+\vai{5}\vai{2}+\vai{6}\vai{1}\Big)
\\&&\hspace{1em}
-\,2\,\Big(\vai{0}\vai{7}(\vai{0}-\vai{7})+\vai{3}\vai{4}(\vai{3}-\vai{4})
+\vai{5}\vai{2}(\vai{5}-\vai{2})+\vai{6}\vai{1}(\vai{6}-\vai{1})\Big)
\\&&\hspace{1em}
+ \,2\,\Big(\big(\vai{3}\vai{5}\vai{6} + \vai{0}\vai{5}\vai{6} +
\vai{0}\vai{3}\vai{6} + \vai{0}\vai{3}\vai{5}\big)
-\big(\vai{4}\vai{2}\vai{1} + \vai{7}\vai{2}\vai{1} +
\vai{7}\vai{4}\vai{1} + \vai{7}\vai{4}\vai{2}\big) \Big)
\vspace{1ex}\\
%
\vf_4 &:=& \Big(\vai{0}^2\vai{7}^2 + \vai{3}^2\vai{4}^2 
           + \vai{5}^2\vai{2}^2 + \vai{6}^2\vai{1}^2\Big)
  \,+\, 4\,\Big(\vai{0}\vai{3}\vai{5}\vai{6}
            \,+\, \vai{7}\vai{4}\vai{2}\vai{1}\Big)\\
&&\hspace{1em}
\,-\, 2\,\Big(\vai{0}\vai{7}\,\vai{3}\vai{4} 
  \,+\, \vai{0}\vai{7}\,\vai{5}\vai{2}
  \,+\, \vai{0}\vai{7}\,\vai{6}\vai{1} 
  \,+\, \vai{3}\vai{4}\,\vai{5}\vai{2}
  \,+\, \vai{3}\vai{4}\,\vai{6}\vai{1} 
  \,+\, \vai{5}\vai{2}\,\vai{6}\vai{1}\Big)\,,
\end{array}
\]
we obtain
\[
\vc_1=\,\vf_1^2\,\vf_4
\hspace{2em}\mbox{and}\hspace{2em}
\vc_1-4\vc_0=\,\vf_3^2\,.
\]
In particular, the two solutions for $\vv$ are
\[
\vv_{1,2}\;=\; \frac{1}{2}\left( 1\pm
\sqrt{\frac{\vc_1-4\vc_0}{\vc_1}}\right)
\;=\;
%
%
\frac{1}{2}\left( 1\pm
\frac{\vf_3(\vai{0},\dots,\vai{7})}
{\vf_1(\vai{0},\dots,\vai{7})\sqrt{\vf_4(\vai{0},\dots,\vai{7})}}\right)\,.
\vspace{-1ex}
\]
\par

{\bf Remarks:}
\kitem
\item[(1)]
Note that whenever $\vv$ solves the equation, then so does
$(1-\vv)$. This symmetry may easily be seen in the original 8 equations
by switching the variables $\ve{\di}$ and $\vq{\di}$.
\item[(2)]
The formulas for $\vf_1, \vf_3$, and $\vf_4$ become very natural if we
recall that $\vai{0}, \vai{3}, \vai{5}, \vai{6}$ correspond to
$\vak{-,-,-}$, $\vak{-,+,+}$, $\vak{+,-,+}$, $\vak{+,+,-}$, respectively.
These variables are those which have an even number of plus signs.\\
This fact may be illustrated by imaging the variables $\vak{\pm,\pm,\pm}$ as
sitting in the corners of a cube. Then, $\vai{0}, \vai{3}, \vai{5}, \vai{6}$
correspond to the vertices of one of the two inscribed regular tetrahedra. 
The remaining $\vai{7}, \vai{4}, \vai{2}, \vai{1}$ are contained in
the opposite corners, respectively.
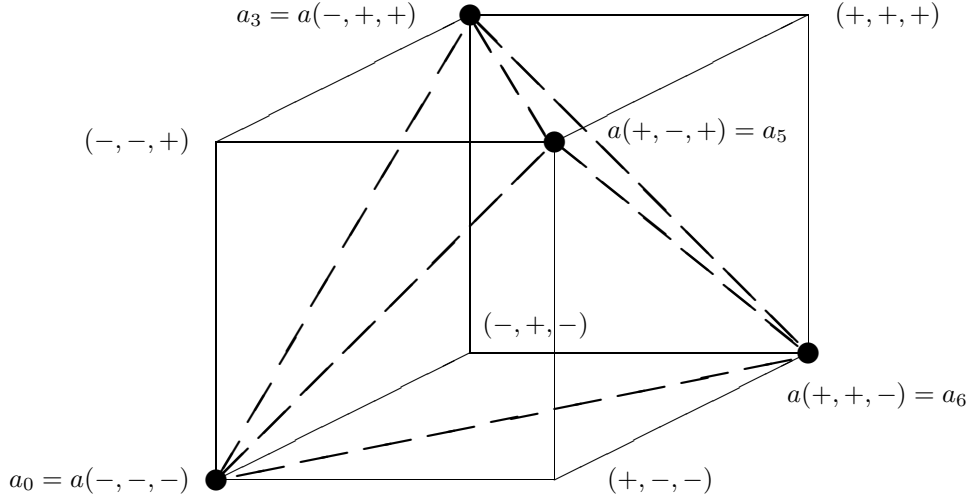
\begin{center}
\unitlength=0.8pt
\begin{picture}(284.00,240.00)(80.00,552.00)
\thinlines
\linethickness{0.1pt}
\put(80.00,560.00){\line(1,0){160.00}}
\put(80.00,720.00){\line(1,0){160.00}}
\put(200.00,780.00){\line(1,0){160.00}}
\put(200.00,620.00){\line(1,0){160.00}}
\put(240.00,560.00){\line(2,1){120.00}}
\put(240.00,720.00){\line(2,1){120.00}}
\put(80.00,720.00){\line(2,1){120.00}}
\put(80.00,560.00){\line(2,1){120.00}}
\put(360.00,620.00){\line(0,1){160.00}}
\put(80.00,560.00){\line(0,1){160.00}}
\put(240.00,560.00){\line(0,1){160.00}}
\put(200.00,620.00){\line(0,1){160.00}}
\thicklines
\multiput(80.00,560.00)(20,20){8}{\line(1,1){15.00}}     
\multiput(90.00,562.00)(30,6.3){9}{\line(5,1){22.00}}    
\multiput(80.00,560.00)(15,27.5){8}{\line(3,5){11.00}}   
\multiput(200.00,780.00)(20,-20){8}{\line(1,-1){15.00}}  
\multiput(202.00,776.00)(12.9,-20){3}{\line(3,-5){8.00}} 
\multiput(240.00,720.00)(20,-16.9){6}{\line(5,-4){15.00}}
\put(80.00,560.00){\circle*{10}}   
\put(240.00,720.00){\circle*{10}}  
\put(200.00,780.00){\circle*{10}}  
\put(360.00,620.00){\circle*{10}}  
%
\put(80.00,560.00){\makebox(0,0)[r]{$\vai{0}=\vak{-,-,-}$\hspace{1em}}}
\put(200.00,780.00){\makebox(0,0)[r]{$\vai{3}=\vak{-,+,+}$\hspace{2em}}}
\put(265.00,725.00){\makebox(0,0)[l]{$\vak{+,-,+}=\vai{5}$}}
\put(350.00,600.00){\makebox(0,0)[l]{$\vak{+,+,-}=\vai{6}$}}
\put(80.00,720.00){\makebox(0,0)[r]{$(-,-,+)$\hspace{1em}}}
\put(240.00,560.00){\makebox(0,0)[l]{\hspace{2em}$(+,-,-)$}}
\put(360.00,780.00){\makebox(0,0)[l]{\hspace{1em}$(+,+,+)$}}
\put(206.00,633.00){\makebox(0,0)[l]{$(-,+,-)$}}
\end{picture}
\end{center}
\item[(3)]
It has been observed by Duco~van~Straten that $\vf_4$ equals the
hyperdeterminant of the three-dimensional matrix $\vA_{\kb\kb\kb}$
formed by the variables $\vak{\pm,\pm,\pm}$, cf.\ Proposition 14.1.7. 
in \cite{GKZ}.
Moreover, $\vf_3$ is a linear combination of the derivatives of $\vf_4$ 
which follows the usual pattern,
\[
2\vf_3=\Big(
\frac{\partial\vf_4}{\partial\vai{0}}+
\frac{\partial\vf_4}{\partial\vai{3}}+
\frac{\partial\vf_4}{\partial\vai{5}}+
\frac{\partial\vf_4}{\partial\vai{6}}\Big) -
\Big(\frac{\partial\vf_4}{\partial\vai{7}}+
\frac{\partial\vf_4}{\partial\vai{4}}+
\frac{\partial\vf_4}{\partial\vai{2}}+
\frac{\partial\vf_4}{\partial\vai{1}}\Big) \,.
\vspace{-2ex}
\]
\kenditem
\par

Finally, we would like to note that the coefficient $\vc_0$ itself does split
into a product of three quadrics:
\[
\vc_0=\vf_{21}\,\vf_{22}\,\vf_{23}
\hspace{1em}\mbox{with}\hspace{1em}
\begin{array}[t]{rcl}
\vf_{21} &:=& (\vai{0}+\vai{4})(\vai{7}+\vai{3}) -
        (\vai{1}+\vai{5})(\vai{6}+\vai{2})\\
\vf_{22} &:=& (\vai{0}+\vai{2})(\vai{7}+\vai{5}) -
        (\vai{4}+\vai{6})(\vai{3}+\vai{1})\\
\vf_{23} &:=& (\vai{0}+\vai{1})(\vai{7}+\vai{6}) -
        (\vai{2}+\vai{3})(\vai{5}+\vai{4})\,.
\vspace{-2ex}
\end{array}
\]
\par


\neu{mmr-elimrest}
Now, we focus on the remaining six variables $\ve{\di}$ and $\vs{\di}$.
Following the above recipe, we obtain again quadratic equations for each of
them, but with much smaller coefficients. They are no longer of degree 6,
but quadratic themselves.
\par

{\bf Notation:}
With $\vA_{\kb\kb\kb}$ being the three-dimensional matrix formed by the
variables $\vak{\pm,\pm,\pm}$, we derive the following ordinary 
$(2\times 2)$ matrices from it:
\kitem
\item
$\vA_+(1):=\vA_{+\kb\kb}$ denotes the layer consisting of the entries
$\vak{+,\kb,\kb}$, i.e., the right hand face of the cube depicted above; 
the remaining (left) one forms the matrix $\vA_-(1):=\vA_{-\kb\kb}$.
Similarly, we may define $\vA_{\pm}(2):=\vA_{\kb\pm\kb}$ and 
$\vA_{\pm}(3):=\vA_{\kb\kb\pm}$.
\item
Considering the sum of the layers, we obtain
$\vA_{\Sigma}(\di):= \vA_{+}(\di) +  \vA_{-}(\di)$ for $\di=1,2,3$.
\kenditem
\par

Using this new terminology, we may recover the quadratic 
$\vc_0$-factors $\vf_{2\di}$
from the end of \zitat{mmr}{elimv} as 
\[
\vf_{2\di} = \det \vA_{\Sigma}(\di)
\hspace{2em}\mbox{with}\hspace{2em}
\di=1,2,3\,.
\]
Fixing a disease index $\di$, the elimination done by {\ktype Singular} tells
us that 
$\ve{\di}$ and $\vq{\di}$ both obey the same quadratic equation. It is
\[
\Big(\det \vA_{\Sigma}(\di)\Big) \, x^2 
- \Big( \det \vA_{\Sigma}(\di) + \det \vA_+(\di) - \det \vA_-(\di) \Big) \, x
+ \Big(\det \vA_+(\di)\Big) = 0\,.
\]
The discriminant is the hyperdeterminant $\,\det \vA=\vf_4$ again. 
Hence, the solutions for $\ve{\di}$ and $\vq{\di}$ are
\[
\Big[(\ve{\di})_{1,2} 
\hspace{0.5em}\mbox{and}\hspace{0.5em} 
(\vq{\di})_{1,2} \Big]
\;=\;
\frac{1}{2}\left( 1+
\frac{\det \vA_+(\di) - \det \vA_-(\di) \pm \sqrt{\det \vA}}
{ \det \vA_{\Sigma}(\di)} \right)
\;=\;
\frac{1}{2}\left(1+\frac{\vg_{2\di}\pm \sqrt{\vf_4}}{\vf_{2\di}}\right)
\]
with $\vg_{2\di}$ being the quadratic polynomials
%
%
\[
\vg_{2\di}\;:=\;
\det\vA_+(\di)-\det\vA_-(\di)=\left\{
\begin{array}{ll}
-\vai{0}\vai{3} + \vai{1}\vai{2} + \vai{4}\vai{7} - \vai{5}\vai{6}
&(\mbox{for}\hspace{0.5em}\di=1)\\
-\vai{0}\vai{5} + \vai{1}\vai{4} + \vai{2}\vai{7} - \vai{3}\vai{6}
&(\mbox{for}\hspace{0.5em}\di=2)\\
-\vai{0}\vai{6} + \vai{1}\vai{7} + \vai{2}\vai{4} - \vai{3}\vai{5}
&(\mbox{for}\hspace{0.5em}\di=3)\,.
\end{array}\right.
\vspace{-1ex}
\]
\par


\neu{mmr-check}
Assuming the general case of $\vf_1\neq 0$, $\,\det \vA     
_{\kb\kb\kb}
\neq 0$, and $\,\det\vA_{\Sigma}(\di)\neq 0$ for each $\di=1,2,3$,
we have narrowed the number of possible values for each of the variables
$\vv, \ve{\di}$, and $\vq{\di}$ down to two.
It remains to check which of the $2^7$
combinations survive to provide an actual solution of the original system
\zitat{mmr}{recall}.\\
This can easily be done by considering the sum of those equations out of 
the original system 
that correspond to a certain face of the cube depicted in \zitat{mmr}{elimv}. 
For instance, adding up
the equations for $\vai{7}, \vai{6}, \vai{5}$, and $\vai{4}$ provides
\[
\vai{7} + \vai{6} + \vai{5} + \vai{4} \;=\; 
\vf_1\,\vv\,\vq{1} \;+\; \vf_1\,(1-\vv)\,\ve{1}\,.
\]
All variables have been eliminated except $\vv$, $\vq{1}$, and $\ve{1}$.
This allows us to show that the $\ve{}$'s must not
equal the $\vq{}$'s. (Assuming $\ve{1}=\vq{1}$, we would obtain 
$\vai{7}+\vai{6}+\vai{5}+\vai{4}=\vf_1\ve{1}$. However, substituting this
value of $\ve{1}$ into the quadratic equation of \zitat{mmr}{elimrest}
yields
\[
\vf_1^2\, \Big(\vf_{21}\,\ve{1}^2-(\vf_{21}+\vg_{21})\,\ve{1}+\det \vA_+(1)\Big)
\;=\;
-\,\vf_{22}\,\vf_{23}\,,
\]
which is generally different from zero.)\\
Now, by Remark \zitat{mmr}{elimv}(1), we may assume that, w.l.o.g.,
$\vv=(\vf_1\sqrt{\vf_4}+\vf_3)/(2\vf_1\sqrt{\vf_4})$. Hence, 
with $\ve{1}=(\vf_{21}+\vg_{21}\mp\sqrt{\vf_4})/(2\vf_{21})$
and $\vq{1}=(\vf_{21}+\vg_{21}\pm\sqrt{\vf_4})/(2\vf_{21})$,
the above equation multiplied with $4\vf_{21}\sqrt{\vf_4}$ becomes
\[
\begin{array}{r@{\hspace{0.5em}}c@{\hspace{0.5em}}l}
\displaystyle
4\,\vf_{21}\,\sqrt{\vf_4}\, \Big(\sum_{\dk=4}^7\vai{\dk} \Big)
&=&
\displaystyle
\Big(\vf_1\sqrt{\vf_4}+\vf_3\Big)
\Big(\vf_{21}+\vg_{21}\pm\sqrt{\vf_4}\Big)
+
\Big(\vf_1\sqrt{\vf_4}-\vf_3\Big)
\Big(\vf_{21}+\vg_{21}\mp\sqrt{\vf_4}\Big)\\
&=&
2\,\vf_1\,\sqrt{\vf_4}\,\big(\vf_{21}+\vg_{21}\big) \,\pm\,
2\,\vf_3\,\sqrt{\vf_4}\,.
\end{array}
\]
In particular, since 
$\,2\vf_{21} \big(\sum_{\dk=4}^7\vai{\dk}\big)
= \vf_1\big(\vf_{21}+\vg_{21}\big)+\vf_3$, 
only the signs on top survive in the formulas of $\ve{1}$ and $\vq{1}$.\\
Finally, one may use {\ktype Singular} again for checking that these values,
together with the similar ones for the remaining variables, 
indeed yield a solution of the original system. This means that we have shown 
the following
\par

{\bf Theorem:}
{\em
If $\vf_1,\vf_4,\vf_{2\di}\neq 0$ for $\di=1,2,3$, then the polynomial system
of \zitat{mmr}{recall}, with the adaption $\vpi{\dk}=\vai{\dk}/\vn$ made in
\zitat{mmr}{prog}, has exactly two solutions. They are
\[
\vv=\frac{\vf_1\,\sqrt{\vf_4}\,\pm\,\vf_3}{2\,\vf_1\,\sqrt{\vf_4}}
\,,\hspace{1em}
\ve{\di}=\frac{\vf_{2\di}+\vg_{2\di}\mp\sqrt{\vf_4}}{2\,\vf_{2\di}}
\,,\hspace{1em}
\vq{\di}=\frac{\vf_{2\di}+\vg_{2\di}\pm\sqrt{\vf_4}}{2\,\vf_{2\di}}
\hspace{2em} (\di=1,2,3).
\]
If some of the above polynomials $\vf_\kb$ do vanish, then the system 
\zitat{mmr}{recall} might have infinitly many solutions or no solution at all.
}
\par

%
%
\section{The MMR coverage}\label{app}


\neu{app-range}
If we apply the previous theory to our statistical problem of estimating 
the MMR coverage, then $\vai{\dk}$ stands for the number of persons of
a prefixed age group observed to have antibody status $\dk$ 
($\dk=0,\dots,7$). Thus, $\vf_1$ is the size of the cohort, and
this number is automatically positive.
On the other hand, we would like to interpret the solutions 
$\vv,\ve{\di},\vq{\di}$, and $\vs{\di}$ of the MMR system as estimations
of the probabilities described in \zitat{mmr}{recall}.
In particular, they should be real numbers and, moreover, 
be contained in the interval $[0,1]$.\\
While in \cite{Gay} the latter is forced by the numerical program used 
to solve the system, our solutions may not have these properties. However,
this should not be considered problematic, but a feature of our method. If the 
solutions fall out of the range making sense, this is a strong hint that
the input data $\vai{\dk}$ are of poor quality.
\par

 
\neu{app-theorem}
In the following, we will formulate the conditions the input data have to
fulfill for yielding apropriate results. Moreover, we will see that,
in the statistical context, only one of the two solutions mentioned in Theorem
\zitat{mmr}{check} survives.
\par

{\bf Theorem:}
{\em
Let $\vai{\dk}$ be the observed number of people in a fixed age group
with antibody status $\dk$. 
Then, the MMR system has a good statistical solution if and only
if
\[
\vf_4(\ku{\va})>0
\hspace{1em}\mbox{and}\hspace{1.2em}
\vf_{2\di}(\ku{\va})\geq \sqrt{\vf_4(\ku{\va})} + 
\big|\vg_{2\di}(\ku{\va})\big|
\hspace{1em}
(\di=1,2,3)\,.
\]
If these conditions are satisfied, then the estimation for 
$\vv,\ve{\di},\vs{\di}$ is
\[
\vv=\frac{\vf_1\,\sqrt{\vf_4}\,+\,\vf_3}{2\,\vf_1\,\sqrt{\vf_4}}
\,,\hspace{1em}
\ve{\di}=\frac{\vf_{2\di}+\vg_{2\di}-\sqrt{\vf_4}}{2\,\vf_{2\di}}
\,,\hspace{1em}
\vs{\di}=\frac{2\,\sqrt{\vf_4}}{\vf_{2\di}-\vg_{2\di}+\sqrt{\vf_4}}
\hspace{2em} (\di=1,2,3).
\vspace{-2ex}
\]
}
\par

{\bf Proof:}
Positivity of $\vf_4$ means that the solutions described in 
Theorem \zitat{mmr}{check} are real. Assuming this, we have
\[
\vv\in [0,1]
\;\Longleftrightarrow\;
\vf_1\sqrt{\vf_4}\pm\vf_3\geq 0
\;\Longleftrightarrow\;
\vf_1^2\vf_4\geq\vf_3^2\,.
\]
On the other hand, we have seen in \zitat{mmr}{elimv} that
\[
\vf_1^2\vf_4 \;=\; \vc_1 \;=\; (\vc_1-4\vc_0)+4\vc_0
\;=\; \vf_3^2\,+\,4\,\vf_{21}\,\vf_{22}\,\vf_{23}\,.
\]
Hence, the condition ``$\,v\in [0,1]\,$'' is equivalent to
$\,\vf_{21}\vf_{22}\vf_{23}>0$.
\vspace{1ex}\\
Since $\,\vs{\di}=(\vq{\di}-\ve{\di})/(1-\ve{\di})$, we know that
\[
\ve{\di},\,\vs{\di}\in [0,1]
\hspace{0.7em}\Longleftrightarrow\hspace{0.7em}
0\,\leq \,\ve{\di}\,\leq \,\vq{\di}\,\leq\, 1\,.
\]
From Theorem \zitat{mmr}{check} we obtain, depending on the choice
of the solution, that
$\,\vq{\di}-\ve{\di} = \sqrt{\vf_4}/\vf_{2\di}$ for $\di=1,2,3\,$
or that $\vq{\di}-\ve{\di} = -\sqrt{\vf_4}/\vf_{2\di}$ for $\di=1,2,3$. 
Anyway, for $\,\vq{\di}\geq\ve{\di}$, 
the polynomials $\vf_{21},\vf_{22},\vf_{23}$ must have the same sign.
Together with $\vf_{21}\vf_{22}\vf_{23}>0$ obtained above,
this means that $\vf_{21},\vf_{22},\vf_{23}>0$. In particular, looking at
Theorem \zitat{mmr}{check}, only the solution with the top sign survives.\\
Finally, it is easy to see that the conditions $\,\ve{\di}\geq 0\,$
and $\,\vq{\di}\leq 1\,$ translate into
$\,\vf_{2\di}\geq \sqrt{\vf_4} - \vg_{2\di}\,$ and
$\,\vf_{2\di}\geq \sqrt{\vf_4} + \vg_{2\di}\,$, respectively.
\hfill$\Box$
\par

 
\neu{app-v}
{\bf Remark:}
If one is only interested in the MMR coverage $\vv$, then the conditions
ensuring a meaningful result may be weakened.  
It follows from the proof of the previous theorem that
\[
\vf_4(\ku{\va})>0
\hspace{1em}\mbox{and}\hspace{1.2em}
\vf_{2\di}(\ku{\va})> 0
\hspace{1em}
(\di=1,2,3)\,.
\]
will do.
\par

%
%
\section{Data}\label{data}

 
\neu{data-germ1}
To illustrate our results, we have chosen some data of some country of the ESEN
Project, \cite{Gay}. These data have not yet been finalized as 
they might be changed according to a new standardization between
the European countries. 
For that reason, the use of these data here is for illustrative purposes
only.\\
The input, i.e., the sampled variables $\vai{\dk}$, may be found 
in the table \zitat{data}{germ2}.
The first table compares our estimation of $\vv, \ve{1}, \ve{2}$, and
$\ve{3}$ by age groups (AG) with that obtained by Gay in \cite{Gay}; 
the variables pointing to his values carry a tilde.
\vspace{-1ex}
\par

%
\begin{center}
\begin{tabular}{|c||r|r||r|r||r|r||r|r||r|r|r|}
\hline
AG&\multicolumn{1}{c|}{$\vvg$}&\multicolumn{1}{c||}{$\vv$}&\multicolumn{1}{c|}{$\veg_1$}&\multicolumn{1}{c||}{$\ve{1}$}&\multicolumn{1}{c|}{$\veg_2$}&\multicolumn{1}{c||}{$\ve{2}$}&\multicolumn{1}{c|}{$\veg_3$}&\multicolumn{1}{c||}{$\ve{3}$}&\multicolumn{1}{c|}{$\vs{1}$}&\multicolumn{1}{c|}{$\vs{2}$}&\multicolumn{1}{c|}{$\vs{3}$}\\
\hline
1&0.227&0.227&0.003&0.005&0.019&0.019&0.014&0.011&0.950&0.861&0.974\\
2&0.642&0.642&0.122&0.144&0.020&0.017&0.090&0.090&0.976&0.878&0.922\\
3&0.715&0.710&0.122&0.112&0.041&0.046&0.090&0.087&1.002&0.912&0.930\\
4&0.837&0.824&0.251&0.279&0.041&0.054&0.106&0.219&1.003&0.886&0.922\\
5&0.859&0.863&0.292&0.252&0.241&0.227&0.106&0.000&1.000&0.886&0.921\\
\hline
6&0.794&0.889&0.621&0.427&0.324&0.094&0.106&-0.037&0.961&0.855&0.830\\
7&0.645&0.847&0.756&0.550&0.502&0.006&0.256&0.258&0.949&0.938&0.678\\
8&0.662&0.794&0.764&0.652&0.502&0.285&0.411&0.356&0.969&0.877&0.798\\
9&0.576&0.900&0.764&0.588&0.665&0.279&0.481&-0.007&0.833&0.857&0.838\\
10&0.478&0.940&0.906&0.667&0.734&0.049&0.631&0.450&0.906&0.892&0.660\\
\hline
\end{tabular}
\end{center}
\par

The main difference between Gay's and our results can be found in the values of 
$\vv, \ve{1}, \ve{2}, \ve{3}$ in the higher age groups.\\ 
Moreover, while Gay has assumed age independent seroconversion rates, 
our solutions $\vs{1}, \vs{2}, \vs{3}$ do vary with age; the most
striking example is the rubella seroconversion $\vs{3}$. The comparison 
of Gay's values with the age average of our solutions for 
$\vs{1}, \vs{2}, \vs{3}$ is as follows:
\vspace{-2ex}
\begin{center}
\begin{tabular}[t]{lrrr}
Seroconversion by N.~Gay: &0.989&0.880&0.910\\
Average of our $\vs{1}, \vs{2}, \vs{3}$: &0.955&0.884&0.847
\end{tabular}
\vspace{-1ex}
\end{center}
\par


\neu{data-germ2}
We can use the equations of \zitat{mmr}{recall} to re-calculate
the expected antibody prevalence out of the solutions obtained for
$\vv, \ve{\di}, \vs{\di}$. In other words, for each antibody status
$(\pm,\pm,\pm)$ we are looking for the number of
people that should have been observed to yield the desired result.\\
Because we used an exact method, it is no surprise that our solutions 
give exactly back the input data; they fill the $\vai{\dk}$-columns in 
the following table. 
On the other hand, using Gay's solutions, we obtain different values
which are contained in the $\vag{\dk}$-columns:
\vspace{-1ex}
\par

\newcommand{\hh}{\hspace{0.2em}}
\begin{center}
\begin{tabular}{|r|r||r|r||r|r||r|r||r|r||r|r||r|r||r|r|}
\hline
\multicolumn{2}{|c||}{$---$}&\multicolumn{2}{c||}{$--+$}&
\multicolumn{2}{c||}{$-+-$}&\multicolumn{2}{c||}{$-++$}&
\multicolumn{2}{c||}{$+--$}&\multicolumn{2}{c||}{$+-+$}&
\multicolumn{2}{c||}{$++-$}&\multicolumn{2}{c|}{$+++$}\\
\hline
$\vag{0}$&$\vai{0}$&$\vag{1}$&$\vai{1}$&
$\vag{2}$&$\vai{2}$&$\vag{3}$&$\vai{3}$&
$\vag{4}$&$\vai{4}$&$\vag{5}$&$\vai{5}$&
$\vag{6}$&$\vai{6}$&$\vag{7}$&$\vai{7}$\\
\hline
\makebox[1.2em][r]{155.8}&\makebox[1.1em][r]{156}&2.3&2&3.1&3&0.5&2&1.0&1&5.0&6&3.7&1&37.7&38\\
49.1&48&5.0&5&1.1&1&1.0&2&7.9&9&12.7&13&8.2&7&90.2&90\\
40.8&42&4.2&4&1.8&2&1.2&0&6.9&6&14.6&11&9.8&8&107.6&114\\
20.1&18&2.5&5&1.0&1&1.2&0&8.2&8&17.7&18&11.6&9&129.7&133\\
14.6&17&1.8&0&4.7&5&1.8&0&7.3&7&16.1&15&15.3&15&153.4&156\\
\hline
10.2&13&1.3&0&5.0&2&1.2&3&17.9&14&14.8&20&20.7&30&145.9&135\\
6.9&11&2.4&4&7.0&1&2.7&3&21.9&16&15.1&13&30.3&40&128.7&127\\
5.0&7&3.5&4&5.0&3&3.8&3&16.5&15&19.1&20&23.2&25&135.9&135\\
3.4&6&3.1&1&6.7&4&6.5&9&11.1&11&14.4&14&26.7&27&122.1&122\\
0.9&2&1.5&2&2.4&1&4.2&4&8.5&7&17.1&17&26.1&28&121.2&121\\
\hline
\end{tabular}
\end{center}
\par


\neu{data-germdisc}
In the following, we will discuss some of the properties of our solutions.
\kitem
\item[(1)]
One should not so much worry about negative rates or rates above 1
as they appear among the $\ve{\di}$ or $\vs{\di}$. In all those cases,
the values are very close to the allowed range.
\item[(2)]
Our major concern is caused by the exposition factors $\ve{2}$ and
$\ve{3}$. They seem to be very small in the higher age groups and,
additionally, they do not increase with age.\\
For the latter, however, we may use the same explanation as Gay did for the
decline of his $\vv$ in older cohorts in that 
the data arise from {\em different}
cohorts in each age group. 
\item[(3)]
As already mentioned before, we did not ad hoc assume that the
seroconversions $\vs{\di}$ are age independent. However, as a result of our
calculations, we obtained values for mumps and measles that did not greatly
vary -- and the averages are quite close to Gay's values.\\
On the other hand, the seroconversion factor for rubella shows an unusual
behavior in the higher age groups and we would be interested in an
explanation for it.
\kenditem
\par

The major difference between Gay's and our approach is the following:
\vspace{0.5ex}\\
{\em Altmann:}
We consider each age group separately; this yields a system of 7 equations
in 7 variables for each group, allowing exact solutions with easy formulas.
\vspace{0.5ex}\\
{\em Gay:}
He considers 10 age groups at once, yielding a system with 70 equations
in 70 variables. Moreover, he creates additional restrictions by
\vspace{-1ex}
\kitem
\item
assuming that the seroconversion $\vs{\di}$ is age independent
(meaning to lose $27$ variables),
\vspace{-1ex}
\item
and by forcing the exposition factors $\ve{\di}$ to increase with age
(meaning to introduce additional inequalities).
\vspace{-2ex}
\kenditem
For the remaining system, Gay uses a numerical approach to find values
for $\vv(\mbox{age})$, $\ve{\di}(\mbox{age})$, and $\vs{\di}$ to fit
into the system as best as possible. Exact solutions are 
of course out of range.
\vspace{1ex}
\par

Thus, the fact that the above problem (2) does not occur in Gay's solutions 
is no surprise at all. It was part of his method to force
all these properties which are, however,
biologically plausible. 
An advantage of Gay's method is that imperfect data in single
age groups might be corrected by the better ones.\\
On the other hand, our method tells which data are better or worse and
gives information about their quality. Moreover, besides exactness,
the main advantage of our approach seems to be that the formulas for
$\vv, \ve{\di}, \vs{\di}$ are mutually independent. Hence, even if one
dislikes the results for the $\ve{\di}$'s or $\vs{\di}$'s, one has still an
explicit formula for the MMR coverage $\vv$ which works well.
\par

%
%

\vspace{3ex}

{\small
\parbox{9cm}{
Doris Altmann\\
Robert Koch Institut\\
Stresemannstr.~90-102\\
D-10963 Berlin, Germany\\
e-mail: altmannd@rki.de}%
\setbox0\hbox{e-mail: altmann@mathematik.hu-berlin.de}\hfill\parbox{\wd0}{
Klaus Altmann\\
Institut f\"ur Reine Mathematik\\
Humboldt-Universit\"at zu Berlin\\
Ziegelstr.~13A\\
D-10099 Berlin, Germany\\
e-mail: altmann@mathematik.hu-berlin.de}}


{\small
\begin{thebibliography}{BCKS}

\bibitem[Ga]{Gay} Gay, N.: A Method for Estimating Coverage of a Multivalent
Vaccine from Antibody Prevalence Data: application to MMR vaccine in 3
European countries. Draft.

\bibitem[GKZ]{GKZ} Gelfand, I.M., Kapranov, M.M., Zelevinsky, A.V.:
Discriminants, Resultants, and Multidimensional Determinants.
Birkh\"auser Boston 1994

\bibitem[GPS]{Singular} Greuel, G.-M., Pfister, G., Sch\"onemann, H.:
Singular. 
System for computer algebra, university of Kaisers\-lautern,
available via {\tt www.mathematik.uni-kl.de}

\end{thebibliography}
}
\end{document}